\documentclass[10pt]{article}
\usepackage[utf8]{inputenc}
\usepackage{latexsym,amssymb,amsmath,url,authblk,skull}

\def\N{\mathbb{N}}

\def\proof{\par\noindent{\em Proof. }}
\def\eproof{\hfill{$\Box$}\bigskip}

\def\ds{\dots}
\def\sus{\subset}
\def\al{\alpha}
\def\be{\beta}
\def\ga{\gamma}

\def\cc{\colon}

\newtheorem{thm}{Theorem}[section]
\newtheorem{prop}[thm]{Proposition}

\newtheorem{lem}[thm]{Lemma}
\newtheorem{prob}[thm]{Problem}
\newtheorem{defi}[thm]{Definition}

\title{Cantor digraphs and abbreviations of formulas}
\author{Martin Klazar }
\date{\today}
\begin{document}
\maketitle
\begin{abstract}  
A~digraph $D=\langle V,E\rangle$ ($E\sus V\times V$) is Cantor if Cantor's theorem\,---\,for no 
set there is a~surjection from it to its power set\,---\,holds in $D$, in 
the sense we explain. We construct a~ZF formula $\varphi$ with 
length $494$ such that $D\models\varphi$ iff $D$ is Cantor. 
In order to obtain $\varphi$, which is a~word over the alphabet 
$$
\{x_1,\,x_2,\,\ds\}\cup
\{\in,\,=,\,\neg,
\,\to,\,\leftrightarrow,\,\wedge,\,
\vee,\,\exists,\,\forall,\,(,\,)\}\,,
$$
we devise abbreviation schemes of ZF formulas. We introduce 
extensive and strongly extensive digraphs and show, by the standard 
argument, that they are Cantor. We construct 
a~countable strongly extensive digraph
with arbitrarily large finite
in-degrees.  
\end{abstract}

\medskip
\begin{quote}
{\footnotesize
Z PLR do MLR jel jsem p\v res \v CSSR,\\
SNB m\'e DKW si stoplo na TK.\footnote{See 
Appendix~\ref{app_zkr} for translation of this Czech sentence.}
\\
$\vdots$

\noindent
Ivan Ml\'adek, song {\em Zkratky} (Acronyms) \cite{zkrat}, 1980s
}
\end{quote}

\section{Introduction\,---\,Cantor digraphs}\label{sec_intro}

In this Section~\ref{sec_intro} we introduce Cantor digraphs and give an overview of the 
article. Section~\ref{sec_ZFform} reviews ZF formulas. In 
Section~\ref{sec_ZFprimeform} we introduce ZF$'$ formulas. This is an 
extension of ZF formulas needed for abbreviation schemes of ZF formulas. 
We develop a~theory of these schemes in Section~\ref{sec_abbrev}. Using an
abbreviation scheme of length~9, in Section~\ref{sec_explCAT} we obtain
a~ZF sentence $\varphi$ with length $494$ such that for every digraph $D$,
$$
\text{$D\models\varphi$ if and only if $D$ is Cantor, i.e. Cantor's theorem holds in $D$}\,.
$$
In Section~\ref{sec_digr} we introduce 
extensive and strongly extensive digraphs and show, by the standard 
argument, that they are Cantor. We construct 
a~countable strongly extensive digraph
with arbitrarily large finite
in-degrees. Section~\ref{sec_concl} contains concluding remarks, also on the relation of this article to G\"odel's completeness theorem.

Classical Cantor's theorem \cite{cant} is a~milestone in set theory.

\begin{thm}[Cantor]\label{thm_Cantor}
For no set $x$ there is a~surjection from $x$ to the 
power set $\mathcal{P}(x)$.    
\end{thm}
(A~function $f\cc A\to B$ is called a~surjection if for every $b\in B$
there is $a\in A$ such that $f(a)=b$.) We embark from the observation 
that one can understand Cantor's theorem more broadly as a~statement 
about digraphs. {\em A~digraph} (directed graph) is a~pair 
$$
D=\langle V,\,E\rangle
$$
of a~nonempty set of {\em vertices} $V$ and a~set of 
{\em arrows} $E\sus V\times V$. For $u,v\in V$ we write $uEv$ iff $\langle u,v\rangle\in E$. One can view $D$ also as a~relational structure 
$$
M=\langle V,\,\in_D,\,=_D\rangle
$$
with {\em the universe} $V$ and two binary 
relations 
$$
\text{$\in_D\,\equiv E$ and 
$=_D\,\equiv\{\langle u,\,u\rangle\cc\;u\in V\}$}\,. 
$$
(The symbol $\equiv$ serves as a~defining equality.) Thus $uEv$ iff $u\in_D v$\,---\,we say that $u$ is {\em a~$D$-element} of $v$\,---\,and the only difference 
between $D$ and $M$ is that for digraphs we understand the equality 
relation $=_D$ implicitly as the usual equality $=$. The satisfaction 
relation $\models$ is understood 
in terms of $M$, 
$$
\text{$D\models\varphi$ really means $M\models\varphi$}\,.
$$
We define $\models$ in detail in the next section.

Let $D=\langle V,E\rangle$ be a~digraph. We introduce
Cantor digraphs via several definitions. For $u\in 
V$ we denote by
$$
N(u)\equiv\{v\in V\cc\;vEu\}
$$
the set of {\em in-neighbors of $u$}; these are the $D$-elements of 
$u$. We define nine digraph predicates SUS, SI, $\ds$, SUR; we use them in  Section~\ref{sec_explCAT}.

\begin{enumerate}
\item For $u,v\in V$ we write $u\sus_D v$ or $\mathrm{SUS}(u;v)$
and say that $u$ is a~{\em $D$-subset} of $v$ if $N(u)\sus N(v)$.

\item For $u,v\in V$ we write $\mathrm{SI}(u;v)$ and say that 
{\em $v$ is the only $D$-element 
of $u$} if $N(u)=\{v\}$.

\item For $u,v\in V$ we write $\mathrm{SIN}(u;v)$ and say that $u$ is {\em the singleton $\{v\}_D$} if $u$ is the only vertex in $V$
such that $N(u)=\{v\}$.

\item For $u,v,w\in V$ we write $\mathrm{DO}(u;v;w)$ and say that
{\em $v$ and $w$ are the only $D$-elements of $u$} if $N(u)=\{v,w\}$. 

\item For $u,v,w\in V$ we write $\mathrm{DOU}(u;v;w)$ and say 
that $u$ is {\em the doubleton 
$\{v,w\}_D$} if $u$ is the only vertex in $V$ such that $N(u)=\{v,w\}$.

\item For $u,v,w\in V$ we write $\mathrm{OPA}(u;v;w)$ and say that $u$ is {\em an ordered pair $\langle v,\,w\rangle_D$} if 
$$
u=\{\{v\}_D,\,\{v,\,w\}_D\}_D\,. 
$$
In Proposition~\ref{prop_OPA} we show that $u=\langle 
v,\,w\rangle_D$ works as the usual ordered pair: $u$ uniquely 
determines $v$ and $w$, 
and vice versa.

\item For $u,v\in V$ we write $\mathrm{REL}(u;v)$ and say that {\em $u$ is a~$D$-relation from $v$ 
to $P(v)$} if every 
$D$-element of $u$ is an ordered pair $\langle w,w'\rangle_D$ 
such that $w\in_D v$ (i.e., $wEv$) and $w'\sus_D v$.

\item For $u,v\in V$ we write $\mathrm{FUN}(u;v)$ and say that {\em $u$ is a~$D$-function from $v$ 
to $P(v)$} if $u$ is a~$D$-relation from $v$ to 
$P(v)$ and for every $D$-element $w$ of $v$ there is exactly 
one ordered pair 
$\langle w,w'\rangle_D$ that is a~$D$-element of $u$.

\item For $u,v\in V$ we write $\mathrm{SUR}(u;v)$ and say that {\em $u$ is a~$D$-surjection from 
$v$ to $P(v)$} if $u$ is a~$D$-function from $v$ to 
$P(v)$ and for every $D$-subset $w'\sus_D v$ there 
is at least one ordered pair $\langle w,w'\rangle_D$ that is 
a~$D$-element of $u$.
\end{enumerate}
We show that definitions 2--6 determine standard ordered pairs. 

\begin{prop}\label{prop_OPA}
Let $D=\langle V,E\rangle$ be a~digraph and let $u,v,w,a,b,c\in V$. Then it is true that
$$
\mathrm{OPA}(u;\,v;\,w)\wedge
\mathrm{OPA}(a;\,b;\,c)\Rightarrow(u=a\iff v=b\wedge w=c)\,.
$$
\end{prop}
\proof
We assume that $u$ is an ordered pair $\langle v,w\rangle_D$ and that
$a$ is an ordered pair $\langle b,c\rangle_D$. Let $\al\equiv u=a$.
Then $N(\al)=\{d,d'\}$ with $N(d)=\{v\}$ and $N(d')=\{v,w\}$.
But also $N(\al)=\{e,e'\}$ with $N(e)=\{b\}$ and $N(e')=\{b,c\}$.
Since $\{d,d'\}=\{e,e'\}$, we get that $v=w$ iff $b=c$, and that $v=b$
and $w=c$.

Let $\be\equiv v=b$ and $\ga\equiv w=c$. Using the uniqueness in 
definitions 3 and 5, we obtain unique
vertices $d$ and $e$ such that $N(d)=\{\be\}$ and $N(e)=
\{\be,\ga\}$. By definition~5, there is a~unique vertex $u=a$ with in-neighbors $\{d,e\}$.
\eproof

The following definition is in fact a~main result of our article. 

\begin{defi}[Cantor digraphs]\label{def_cantDigr}
Let $D=\langle V,E\rangle$ be a~digraph. We say that the digraph 
$D$ is Cantor if for no vertex $u\in V$ there exists a~$D$-surjection 
$v\in V$ from $u$ to $P(u)$. That is,
$$
D\models(\forall x_1\neg(\exists x_2\,\mathrm{SUR}(x_2;\,x_1)))\,.
$$
We also say that Cantor's theorem holds in $D$.
\end{defi}
Not every digraph is Cantor. For example, in the digraph
$$
D=\langle V,\,=_D\rangle\,
$$
where the only arrows are loops, one at each vertex, every vertex $u\in 
V$ is a~surjection from $u$ to $P(u)$. However, by removing 
one or more of these loops, we get a~digraph
in which Cantor's theorem holds.

Let $D=\langle V,E\rangle$ be a~digraph. For any vertex 
$u\in V$ we define  {\em the $D$-power set} of $u$ as
$$
\mathcal{P}(u)\equiv\{v\in V\cc\;v\sus_D u\}=
\{v\in V\cc\;N(v)\sus N(u)\}\ \ (\sus V)\,,
$$
and for any ordered pair $\langle u,v\rangle_D$ ($\in V$) we consider
the corresponding (real) ordered pair $\langle u,v\rangle$
($\in V\times V$).

\begin{prop}
Let $D=\langle V,E\rangle$ be a~digraph, $u,v\in V$ and let $\mathrm{SUR}(u;v)$. Then
$$
\{\langle a,\,b\rangle\cc\;\langle a,\,b\rangle_D\in N(u)\}\ \ (\sus V\times V) 
$$
is a~surjection from $N(v)$ to $\mathcal{P}(v)$.
\end{prop}
\proof
This follows at once from definition~9 of SUR.
\eproof

Thus we do not view the statement of Theorem~\ref{thm_Cantor} informally 
in terms of naive set theory, but we view it formally in terms of digraphs. In this perspective Theorem~\ref{thm_Cantor} 
claims that any digraph 
$$
\text{$D_{\mathrm{ZF}}=\langle V,\,E\rangle$ (i.e., structure
$M_{\mathrm{ZF}}=\langle V,\,\in_D,\,=_D\rangle$)}
$$
with the property that it is a~model of the ZF (Zermelo--Fraenkel) set theory, 
is Cantor. By G\"odel's second incompleteness theorem (\cite[Chapter~IV]{daws}), the existence of
$D_{\mathrm{ZF}}$ cannot be established by formal means inside 
ZF.

The primary theme of this article is to get a~completely rigorous definition of Cantor digraphs by expanding the 
displayed formula in Definition~\ref{def_cantDigr} in a~ZF formula $\varphi$. We accomplish it in 
Sections~\ref{sec_ZFform}--\ref{sec_explCAT}. The formula $\varphi$
is obtained in Theorem~\ref{thm_CantorForm} and is stated explicitly
at the end of Section~\ref{sec_explCAT}. Once we rigorously define
Cantor digraphs, it is a~natural idea to have some examples of them. This
is the secondary, in this article somewhat neglected, theme that is 
treated in Section~\ref{sec_digr}. Here we just mention that
any digraph $D_{\mathrm{ZF}}$ is Cantor due to the ZF axiom schema of specification.

\section[ZF formulas]{ZF formulas} \label{sec_ZFform}

To get the formula $\varphi$ we need a~good grasp of ZF formulas. Let $\N$ be the set $\{1,2,\ds\}$ of natural numbers. 
Recall that {\em a~word} $u$ over an alphabet $A$, which may be any nonempty set, is a~finite sequence
$$
u=\langle a_1,\,a_2,\,\ds,\,a_n\rangle=a_1\,a_2\,\ds\,a_n
$$
($n\in\N$) of elements $a_i$ in $A$,
or {\em the empty word $\lambda$}. The length $n$ of $u$ is denoted by 
$|u|$, and $|\lambda|=0$. If $a$ is in $A$, then $|u|_a$ is the number of occurrences 
of $a$ in $u$, that is, the number of indices $i$ such that $a_i$ is 
$a$. We denote the set of words over $A$ by $A^*$. We shall work with the alphabet
$$
\mathcal{A}\equiv\{x_i\cc\;i\in\N\}\cup
\{\,\in,\,=,\,\neg,
\,\to,\,\leftrightarrow,\,\wedge,\,
\vee,\,\exists,\,\forall,\,(,\,)\,\}\,.
$$ 
It consists of countably many {\em set variables $x_i$} and of 
eleven symbols with well-known meanings.

\begin{defi}\label{def_atZFform}
Atomic {\em ZF} formulas are the words over $\mathcal{A}$ with length $5$ 
$$
\text{$(x_i\in x_j)$ and $(x_i=x_j)$},\ i,\,j\in\N\,.
$$
\end{defi}
\vspace{-3mm}
\begin{defi}\label{def_ZFform}
A~word $u\in\mathcal{A}^*$ is a~{\em ZF} formula if and only if there exists a~finite sequence
$$
u_1,\,u_2,\,\ds,\,u_n
$$
of words $u_i\in\mathcal{A}^*$ such that $u_n$ is $u$ and for every index $i=1,2,\ds,n$ one of eight cases occurs.
\begin{enumerate}
\item The word $u_i$ is an atomic {\em ZF} formula. 
\item There exists an index $j<i$ such that the word $u_i$ has form
$\neg u_j$ and length $1+|u_j|$.
\item There exist indices $j,j'<i$ such that the word $u_i$ has form
$(u_j\to u_{j'})$
and length $3+|u_j|+|u_{j'}|$.
\item There exist indices $j,j'<i$ such that the word $u_i$ has form
$(u_j\leftrightarrow u_{j'})$
and length $3+|u_j|+|u_{j'}|$.
\item There exist indices $j,j'<i$ such that the word $u_i$ has form
$(u_j\wedge u_{j'})$
and length $3+|u_j|+|u_{j'}|$.
\item There exist indices $j,j'<i$ such that the word $u_i$ has form
$(u_j\vee u_{j'})$
and length $3+|u_j|+|u_{j'}|$.
\item There exists an index $j<i$ and an index $k$ such that the word $u_i$ has form
$(\exists x_k u_j)$
and length $4+|u_j|$.
\item There exists an index $j<i$ and an index $k$ such that the word $u_i$ has form
$(\forall x_k u_j)$
and length $4+|u_j|$.
\end{enumerate}
\end{defi}
The sequence $u_1$, $u_2$, $\ds$, $u_n$ is sometimes called 
{\em a~generating word of $u$}. It follows that every word $u_i$ in it is a~ZF  
formula. It is not hard to see that the shortest generating word of $u$ has the property that 
$$
\text{$u_i\ne u_j$ for $i\ne j$ and
every word $u_i$ is a~(contiguous) subword of $u$}\,.
$$
Using this we could easily devise an algorithm that for every input word 
over $\mathcal{A}$ decides if it is a~ZF formula. This is not so clear 
for some other definitions of formulas appearing in the literature. If 
$$
u\equiv a_1\,a_2\,\ds\,a_n
$$ 
is a~ZF formula and $1\le i\le j\le n$, the 
subword $a_ia_{i+1}\ds a_j$ is {\em a~subformula} of $u$ if the 
word
$$
b_1\,b_2\,\ds\,b_{j-i+1}\,,
$$
where $b_1\equiv a_i$, $b_2\equiv a_{i+1}$, $\ds$, $b_{j-i+1}\equiv 
a_j$ is a~ZF formula. 

Why do we not shorten Definition~\ref{def_ZFform} by selecting 
a~subset of connectives and quantifiers and then expressing the rest in terms of the selected symbols? For example, 
\cite{shoe} selects $\neg$, $\vee$ and $\exists$, \cite{schi} uses all 
connectives and quantifiers and \cite{kune} selects $\neg$, 
$\wedge$ and $\exists$. In an early version of our article we selected 
$\neg$, $\to$ and $\forall$. However, this minimalism is disadvantageous. It makes the process of abbreviation 
unnecessarily complicated and makes the sought-for formula $\varphi$ 
unnecessarily long.

An aspect of formulas and similar objects like terms, which is 
sometimes neglected, is unique reading lemmas, or URL.

\begin{prop}[URL for ZF formulas]\label{prop_URLform}
Suppose that $u\in\mathcal{A}^*$ is a~{\em ZF} formula. Then exactly one of eight cases occurs.
\begin{enumerate}
\item There is a~unique atomic {\em ZF} formula $v$ such that $u$ is $v$.
\item There is a~unique {\em ZF} formula $v$ such that $u$ is $\neg v$.
\item There are unique {\em ZF} formulas $v$ and $v'$ such that 
$u$ is $(v\to v')$.
\item There are unique {\em ZF} formulas $v$ and $v'$ such that 
$u$ is $(v\leftrightarrow v')$.
\item There are unique {\em ZF} formulas $v$ and $v'$ such that 
$u$ is $(v\wedge v')$.
\item There are unique {\em ZF} formulas $v$ and $v'$ such that 
$u$ is $(v\vee v')$.
\item  There is a~unique {\em ZF} formula $v$ and a~unique index $k$ such that $u$ is $(\exists x_k v)$.
\item  There is a~unique {\em ZF} formula $v$ and a~unique index $k$ such that $u$ is $(\forall x_k v)$.
\end{enumerate}
\end{prop}
Nontrivial cases are the binary ones, 3--6. Sometimes it is suggested that URL 
for formulas and similar objects are automatic corollaries of 
inductive definitions, but this is a~fallacy. In reality URL like
Proposition~\ref{prop_URLform} follow from the next result.

Let 
$$
u=a_1\,a_2\,\ds\,a_n
$$ be a~word over the two-element alphabet 
$\{),(\}$. {\em A~good bracketing of $u$} is a~partition $P$ of 
$\{1,2,\ds,n\}$
in two-element blocks $B=\{i_B<j_B\}$ such that for every $B$ in $P$,
$$
\text{$a_{i_B}$ is $($ and $a_{j_B}$ is $)$}\,,
$$
and that no two blocks $B$ and $C$ in $P$ {\em cross},
$$
\text{neither $i_B<i_C<j_B<j_C$ nor $i_C<i_B<j_C<j_B$}\,.
$$

\begin{prop}\label{prop_brack}
Every word in $\{),(\}^*$ has at most one good bracketing.    
\end{prop}
Alternatively, one can avoid brackets and still have URL by using prefix (Polish) notation, as in 
\cite{shoe}. We leave proofs of Propositions~\ref{prop_URLform} and \ref{prop_brack} as exercises for the interested reader.

Some results on formulas require URL and some can be proven just by induction along generating words. The correct definition of the satisfaction relation $\models$ 
belongs to the former results. The fact that two subformulas of 
a~formula are either disjoint or one contains the other, in particular 
atomic subformulas are disjoint, belongs to the latter results. 

So let us (correctly) define the relation $\models$, in fact relations $\models_f$. We suppose that $D=\langle V,E\rangle$ is a~digraph, $u\in\mathcal{A}^*$ is a~ZF formula and that 
$$
f\cc\{x_i\cc\;i\in\N\}\to V
$$ 
is a~realization of variables by vertices. We proceed by induction on $|u|$
and distinguish eight cases according to Proposition~\ref{prop_URLform}.

\begin{enumerate}
\item If $u$ is $(x_i\in x_j)$ then $D\models_f u$ iff $f(x_i)Ef(x_j)$.
If $u$ is $(x_i=x_j)$ then $D\models_f u$ iff $f(x_i)$ equals $f(x_j)$.
\item If $u$ is $\neg v$ then $D\models_f u$ iff it is not true that $D\models_f v$.
\item If $u$ is $(v\to v')$ then it is not true that $D\models_f u$ iff
it is true that $D\models_f v$ but not that $D\models_f v'$.
\item If $u$ is $(v\leftrightarrow v')$ then 
$D\models_f u$ iff
both $D\models_f v$ and $D\models_f v'$ are, or are not, valid.
\item If $u$ is $(v\wedge v')$ then 
$D\models_f u$ iff
both $D\models_f v$ and $D\models_f v'$ are valid.
\item If $u$ is $(v\vee v')$ then 
$D\models_f u$ iff
at least one of $D\models_f v$ and $D\models_f v'$ is valid.
\item If $u$ is $(\exists x_k v)$ then $D\models_f u$ iff there exists a~map $g\cc\{x_i\cc\;i\in\N\}\to V$
identical with $f$ except (possibly)
for the value $g(x_k)$ such that
$D\models_g v$.
\item If $u$ is $(\forall x_k v)$ then $D\models_f u$ iff for every 
map $g$ as in the previous case
we have $D\models_g v$.
\end{enumerate}
In the next section we define a~family of ZF formulas $u$ 
called sentences for which the validity of $D\models_f u$ does not depend on $f$ and one can write just
$D\models u$.

\section[ZF$'$ formulas]{ZF$'$ formulas} \label{sec_ZFprimeform}

Let $\mathcal{Q}$ be a~nonempty finite set of  {\em predicates} $q$, each of which has an {\em arity} 
$a(q)\in\N$.
Let
$$
\mathcal{V}=
\{\,x,\,y,\,z,\,a,\,b,\,c,\,y_1,
\,y_2,\,\ds\,\}
$$ 
be a~countable set of {\em new set variables} which serve as arguments of predicates. 
We extend the alphabet $\mathcal{A}$ to
$$
\mathcal{A}'\equiv\mathcal{A}
\cup
\{\,;\,\}\cup
\mathcal{V}
\cup\mathcal{Q}\,,
$$
where $;$ is a~symbol for separating arguments of predicates.

\begin{defi}\label{def_ZFprimeAtFor}
Atomic {\em ZF}$'$ formulas are all words over $\mathcal{A}'$ with length $5$ and form
$$
\text{$(\al\in\al')$ and $(\al=\al')$}\,,
$$
where $\al$ and $\al'$ are variables in $\{x_1,x_2,\ds\}\cup\mathcal{V}$, and all words over $\mathcal{A}'$ with length $2a(q)+2$ and form
$$
q(\be_1;\,\be_2;\,\ds;\,\be_{a(q)})\,,
$$
where $q\in\mathcal{Q}$ and $\be_i$ are variables in 
$\{x_1,x_2,\ds\}\cup\mathcal{V}$.
\end{defi}

\begin{defi}
\label{def_ZFprimeFor} 
{\em ZF}$'$ formulas are the words over 
$\mathcal{A}'$ obtained according to the modified Definition~\ref{def_ZFform}. 
We extend case~1 to atomic {\em ZF}$'$ formulas and keep the rest of Definition~\ref{def_ZFform} the same.
\end{defi}
Thus quantification of new variables is not allowed. Subformulas of ZF$'$ formulas
are defined as for ZF formulas.

We review free and bound (occurrences of) variables.
Let $u$ be a~ZF$'$ formula and $\al$ be a variable in 
$\{x_1,x_2,\ds\}\cup\mathcal{V}$. An occurrence of $\al$ in $u$ is {\em 
bound} if it lies in a~subformula of 
$u$ of the quantified form $(\forall\al v)$ or $(\exists\al v)$. Else the 
occurrence of $\al$ is {\em free}. By Definition~\ref{def_ZFprimeFor}, new 
variables have only free occurrences. A~ZF$'$ formula is {\em 
a~sentence} if it has no free occurrence of any variable. It 
follows that every sentence is a~ZF formula. It is not hard to see that 
for ZF sentences $u$ the validity of the satisfaction
$$
D\models_f u
$$
is independent of the realization of variables $f$. We therefore write just
$D\models u$.

\section{Well formed abbreviation schemes}\label{sec_abbrev}

In mathematical logic and set theory, abbreviations of formulas are not treated sufficiently 
rigorously, despite the fact that 
they (should) constitute a~fundamental and indispensable 
syntactic tool. Now we fix it.

We introduce well formed abbreviation schemes and begin with
shortcuts which define predicates.

\begin{defi}\label{def_shortc}
A~shortcut is an expression $\Phi$ of the form
$$
q(y_1;\,y_2;\,\ds;\,y_{a(q)})\equiv
\varphi\,.
$$    
Here $q\in\mathcal{Q}$ and $\varphi$ is a~{\em ZF}$'$ 
formula such that no variable $x_i$
has a~free occurrence in it
and no variable in 
$$
\mathcal{V}\setminus\{y_1,\,y_2,\,\ds,\,y_{a(q)}\}
$$ 
is used. If $a(q)\le3$\,---\,in the next section this will be always the case\,---\,we may use in $q(\ds)$ 
instead of the $y_i$'s the new variables $x$, $y$ and $z$.
\end{defi} 

Suppose that  the predicates in $\mathcal{Q}$ are labeled as 
$$
\{q_1,\,q_2,\,\ds,\,q_l\}\,.
$$
Let $\Phi$ be a~shortcut of the form $q_i(\ds)\equiv\varphi$. We denote by $R(\Phi)$ the set of indices of predicates 
appearing in $\varphi$, and by $V(\Phi)$ the set of indices~$j$ of the 
variables $x_j$ appearing in $\varphi$. If $A,B\sus\N$, we write 
$A<B$ if $m<n$ for every $m$ in $A$ and every $n$ in $B$. In particular, $A<B$ holds if $A$ or $B$ is empty.

\begin{defi}\label{def_wellForAS}
Let the predicates in $\mathcal{Q}$ be labeled as above. A~well formed abbreviation scheme is an $l$-tuple
$$
U\equiv\langle\Phi_1,\,\Phi_2,\,\ds,\,
\Phi_l\rangle
$$
of shortcuts  $\Phi_i$ of the form $q_i(\ds)\equiv\varphi_i$ such that the inequalities 
$$
R(\Phi_i)<\{i\}\,\text{ and }\,
V(\Phi_1)<V(\Phi_2)<\ds<V(\Phi_l)
$$
hold.
\end{defi}
In fact, it suffices to require that the sets $V(\Phi_i)$ are disjoint. The former 
condition $R(\Phi_i)<\{i\}$ is a~natural one, any predicate can 
be defined only in terms of already defined predicates. In other words, definitions of 
predicates must not be circular. 
The latter disjointness condition is a~standard substitutability 
condition, used often for terms. Abbreviation schemes are 
a~more precise and purely syntactic 
version of towers of conservative extensions of theories by definitions of new predicates, as described in the theorem on definition 
of a~predicate in Sochor \cite[str. 21--22]{soch}. 

To define expansion along a~well formed abbreviation scheme, we need 
substitution operations on words. For two natural numbers $m\le n$ we define sets
$$
[n]\equiv\{1,\,2,\,\ds,\,n\}\,\text{ and }\,
[m,\,n]\equiv\{m,\,m+1,\,\ds,\,n\}\,.
$$
We set $[0]\equiv\emptyset$ and $[m,n]\equiv\emptyset$ if $m>n$. 

\begin{defi}\label{def_replSUbw}
Let
$$
u=a_1\,a_2\,\ds\,a_n\,\text{ and }\,v=b_1\,b_2\,\ds\,b_{n'}
$$
be two nonempty words over an alphabet $A$, and let $l$ and $m$
be natural numbers such that $1\le l\le m\le n$. We define the word 
$$
\mathrm{rep}(u,\,v,\,l,\,m)=c_1\,c_2\,\ds\,c_{n''}\ \ (\in A^*)
$$
with length $n''=n-(m-l+1)+n'$ by setting 
\begin{itemize}
\item $c_i\equiv a_i$ for $i\in[l-1]$, 
\item $c_i\equiv b_{i-l+1}$ for $i\in[l,\,l+n'-1]$, and
\item $c_i\equiv a_{i-(l+n')+m+1}$ for $i\in[l+n',n'']$.
\end{itemize}
\end{defi}
Thus one replaces in $u$ the subword at $[l,m]$ with the word $v$. We define two related operations.

\begin{defi}\label{def_symForSym}
Let $u\in A^*\setminus\{\lambda\}$ and $a_i,b_i\in A$ for $i\in[k]$, $k\ge1$, be such that $a_i\ne a_j$ for $i\ne j$. Then
$$
\mathrm{sub}_1(u,\,a_1/b_1,\,
\ds,\,a_k/b_k)
$$
is the word in $A^*$ obtained from $u$
by means of the operation in Definition~\ref{def_replSUbw} 
by replacing for $i\in[k]$ every occurrence of $a_i$ in $u$ with $b_i$.
\end{defi}

For the second operation 
we need a~more detailed version of the operation $\mathrm{rep}(\ds)$.
Let $u$, $v$, $l$ and $m$ be as in Definition~\ref{def_replSUbw}. 
Let $l',m'\in\N$ be such that $1\le l'\le m'\le |u|$ and that the 
intervals $[l,m]$ and $[l',m']$ are disjoint. We define
$$
\text{$\mathrm{rep}_0(u,\,v,\,l,\,m,\,l',\,m')$ to be the triple
$\langle w,\,l'',\,m''\rangle$
}
$$
such that $w\equiv\mathrm{rep}(u,v,l,m)$, $l''\equiv l'$ and $m''\equiv m'$ if 
$m'<l$, and $l''\equiv l'-(m-l+1)+|v|$ and $m''\equiv m'-(m-l+1)+|v|$ if $l'>m$. Thus
we record the action of the replacement on the pair $l',m'$.

\begin{defi}\label{def_symuSubs}
Let $u$, $v_1$, $\ds$, $v_k$, $k\in\N$, be $k+1$ nonempty words over an alphabet $A$ and $l_i,m_i\in\N$ for $i\in[k]$ be such that $1\le 
l_i\le m_i\le |u|$ and that the intervals $[l_i,m_i]$ are pairwise 
disjoint. We consider a~sequence of $k+1$ $(2k+1)$-tuples
$$
\text{$\langle u_i,\,l_{1,i},\,m_{1,i},\,\ds,\,l_{k,i},\,m_{k,i}\rangle$ for $i\in[k+1]$}\,,
$$
starting for $i=1$ with $u_1\equiv u$, $l_{j,1}\equiv l_j$, $m_{j,1}\equiv m_j$, 
and for $i\in[2,k+1]$ and $j\in[i,k]$ continuing with
$$
\langle u_i,\,l_{j,i},\,m_{j,i}
\rangle\equiv\mathrm{rep}_0(u_{i-1},\,v_{i-1},\,l_{i-1,i-1},\,m_{i-1,i-1},\,l_{j,i-1},\,m_{j,i-1})\,.
$$
Then we define
$$
\text{
$\mathrm{sub}_2(u,\,v_1,\,\ds,\,v_k,\,l_1,\,m_1,\,\ds,\,l_k,\,m_k)$ to be the word $u_{k+1}$}\,.
$$
\end{defi}
Thus we replace in the order $i=1,2,\ds,k$ the subword of $u$ at $[l_i,m_i]$ with the word $v_i$. Since the intervals $[l_i,m_i]$ are disjoint, after any 
permutation of the triples
$$
\langle v_1,\,l_1,\,m_1\rangle,\,
\ds,\,\langle v_k,\,l_k,\,m_k\rangle
$$
the operation $\mathrm{sub}_2(\ds)$
yields the same result.

We proceed to expansions along abbreviation schemes. Let 
$$
U=\langle\Phi_1,\,\Phi_2,\,\ds,\,
\Phi_l\rangle
$$
be a~well formed abbreviation scheme in which the $i$-th shortcut $\Phi_i$ is
$$
q_i(y_1;\,y_2;\,\ds;\,y_{a(q_i)})\equiv\varphi_i
$$
(or the arguments in $q_i(\ds)$ are some of $x$, $y$ and $z$). We define 
by induction on $i\in[l]$ the expansion $E_i$ of $\Phi_i$ along 
$U$. It is a~unique ZF$'$ formula free of predicates. For $i=1$ we set 
$E_1\equiv\varphi_1$. Since $U$ is well formed, $R(\Phi_1)=\emptyset$ 
and $E_1$ is indeed free of predicates. 

We suppose that $i>1$ is in $[l]$ and that the expansions $E_1$, 
$E_2$, $\ds$, $E_{i-1}$ are already defined (they are ZF$'$ formulas free of predicates). To get $E_i$, we 
find all atomic subformulas of 
$\varphi_i$ involving a~predicate. They are determined by the triples
$$
\langle k_1,\,l_1,\,m_1\rangle,\,
\langle k_2,\,l_2,\,m_2\rangle,\,\ds,\,
\langle k_s,\,l_s,\,m_s\rangle
$$
such that $s\ge0$, $1\le l_j\le m_j\le|\varphi_i|$ and 
the subword of $\varphi_i$ at $[l_j,m_j]$ is an atomic subformula 
involving $q_{k_j}$. From the remark in Section~\ref{sec_ZFform} we know that the intervals 
$[l_j,m_j]$ are mutually disjoint. 
If $s=0$ ($\varphi_i$ contains no predicate), we set $E_i$ to be 
$\varphi_i$ and are done. If $s\ge1$, then $k_j<i$ for every $j\in[s]$
because $U$ is well formed.

Let $s\ge 1$, $j$ run in $[s]$ and let the subword of $\varphi_i$ at $[l_j,m_j]$ be
$$
q_{k_j}(\al_1\,;\al_2;\,\ds;\,\al_{a(q_{k_j})})
$$
where the $\al_t$ are variables in 
$\{x_1,x_2,\ds\}\cup\mathcal{V}$. Suppose that the left-hand side of $\Phi_{k_j}$ ($k_j<i$) is
$$
q_{k_j}(\be_1\,;\be_2;\,\ds;\,\be_{a(q_{k_j})})
$$
where the $\be_t$ are variables in 
$\mathcal{V}$. Using Definition~\ref{def_symForSym} we set
$$
M_j\equiv\mathrm{sub}_1(E_{k_j},\,
\be_1/\al_1,\,\ds,\,
\be_{a(q_{k_j})}/\al_{a(q_{k_j})})
$$
and using Definition~\ref{def_symuSubs} we set
$$
E_i\equiv\mathrm{sub}_2(\varphi_i,\,M_1,\,\ds,\,M_s,\,l_1,\,m_1,\,\ds,\,
l_s,m_s)\,.
$$

\begin{defi}\label{def_expaAlon}
The {\em ZF}$'$ formulas 
$$
E_1,\ E_2,\ \ds,\ E_l 
$$ 
obtained are free of predicates and we call them 
expansions (along the well 
formed abbreviation scheme $U$).     
\end{defi}
The result of the expansion process is typically the last expansion 
$E_l$. We call expansion
in the order $E_1$, $E_2$, $\ds$, $E_l$ the {\em forward expansion}.
We exemplify it in  the next section.
One could define expansion also
in the opposite order, starting from $E_l$ and going to the $E_i$ 
with $i<l$, but we do not consider this possibility here. 

\begin{prop}\label{prop_freeVar}
Let $i\in[l]$. We characterize occurrences of variables in
expansions $E_i$. Every occurrence of 
every variable $x_j$ is bound. The only free occurrences are  
of (some of) the variables in $\mathcal{V}$ used in the left-hand 
side $q_i(\ds)$ of $\Phi_i$.
\end{prop}
\proof
This follows from the definition of shortcuts in Definition~\ref{def_shortc} and from the expansion process.
\eproof

\section{The sentence $\varphi$}\label{sec_explCAT}

We obtain a~ZF sentence $\varphi$ such that for every digraph $D$,
$$
D\models\varphi\iff\text{$D$ is Cantor}\,.
$$
This sentence arises by expanding a~well formed abbreviation scheme 
$$
U_0\equiv\langle\Phi_1,\,\Phi_2,\,\ds,\,\Phi_9\rangle
$$ 
described below. It uses predicates
$$
\mathcal{Q}\equiv\{q_1,\,q_2,\,\ds,\,q_9\}\equiv
\{\mathrm{SUS},\,\mathrm{SI},\,\mathrm{SIN},\,
\mathrm{DO},\,\mathrm{DOU},\,\mathrm{OPA},\,\mathrm{REL},\,\mathrm{FUN},\,\mathrm{SUR}\}
$$
(respectively), defined already in Section~\ref{sec_intro}. Let 
$D=\langle V,E\rangle$ be any digraph and 
$$
f\cc\{x_1,\,x_2,\,\ds\}\cup\mathcal{V}\to V
$$
be any realization of variables by vertices. Besides the length of 
expansions, we keep (just of interest) track of the number of negations 
used.

\begin{lem}[$\Phi_1$]\label{lem_SUS}
The shortcut $\Phi_1$ is
$$
q_1(x;\,y)\equiv\mathrm{SUS}(x;\,y)\equiv
(\forall x_1((x_1\in x)\to(x_1\in y)))\,.
$$
Thus $|E_1|=17$, $|E_1|_{\neg}=0$ and
$$
\text{$D\models_f\mathrm{SUS}(x;\,y)$ if and only if $f(x)\sus_D f(y)$}\,.
$$
\end{lem}
\proof
The syntactic part is an easy count $4+3+5+5=17$ and the fact that  
no $\neg$ was used. Semantically, the satisfaction in $D$ of the formula  $\mathrm{SUS}(x,y)$ matches the word description. 
\eproof 

\noindent
Again,
{\footnotesize
$$
E_1\equiv
(\forall x_1((x_1\in x)\to(x_1\in y)))\,.
$$
}

\begin{lem}[$\Phi_2$]\label{lem_SI}
The shortcut $\Phi_2$ is
$$
q_2(x;\,y)\equiv\mathrm{SI}(x;\,y)\equiv
(\forall x_2((x_2\in x)\leftrightarrow(x_2=y)))\,.  
$$
Thus $|E_2|=17$, $|E_2|_{\neg}=0$ and 
$$
\text{$D\models_f\mathrm{SI}(x;\,y)$ if and only if $f(y)$ is the only $D$-element of $f(x)$}\,.
$$ 
\end{lem}
\proof
The syntactic part is an easy 
count $4+3+5+5=17$ and observation that no negation was used. Semantically, the satisfaction in $D$ of the formula 
$\mathrm{SI}(x;y)$ matches the word description.
\eproof

\noindent
Again,
{\footnotesize
$$
E_2\equiv
(\forall x_2((x_2\in x)\leftrightarrow(x_2=y)))\,.
$$
}

\begin{lem}[$\Phi_3$]\label{lem_SIN}
The shortcut $\Phi_3$ is
$$
q_3(x;\,y)\equiv\mathrm{SIN}(x;\,y)\equiv
(\forall x_3(\mathrm{SI}(x_3;\,y)\leftrightarrow(x_3=x)))\,.
$$
Thus $|E_3|=29$, $|E_3|_{\neg}=0$ and 
$$
\text{$D\models_f\mathrm{SIN}(x;\,y)$ if and only if $f(x)$ is the singleton $\{f(y)\}_D$}\,.
$$
\end{lem}
\proof
In view of Lemma~\ref{lem_SI}, 
the syntactic part is an easy count $4+3+17+5=29$ and observation that no negation was used. Semantically, the 
satisfaction in $D$ of the formula $\mathrm{SIN}(x;y)$ matches the word description.
\eproof

\noindent
By the forward expansion, 
{\footnotesize
$$
E_3\equiv
(\forall x_3((\forall x_2((x_2\in x_3)\leftrightarrow(x_2=y)))\leftrightarrow(x_3=x)))\,.
$$
}

\begin{lem}[$\Phi_4$]\label{lem_DO}
The shortcut $\Phi_4$ is
$$
q_4(x;\,y;\,z)\equiv\mathrm{DO}(x;\,y;\,z)\equiv
(\forall x_4((x_4\in x)\leftrightarrow((x_4=y)\vee(x_4=z))))\,.
$$
Thus $|E_4|=25$, $|E_4|_{\neg}=0$ and
$$
\text{$D\models_f\mathrm{DO}(x;\,y;\,z)$ iff $f(y)$ and $f(z)$ are the only $D$-elements of $f(x)$}\,.
$$
\end{lem}
\proof
The syntactic part is an easy 
count $4+3+5+3+5+5=25$ and observation that no negation was 
used. Semantically, the satisfaction in $D$ of the 
formula $\mathrm{DO}(x;y;z)$ matches the word description.
\eproof

\noindent
Again,
{\footnotesize
$$
E_4\equiv
(\forall x_4((x_4\in x)\leftrightarrow((x_4=y)\vee(x_4=z))))\,.
$$
}

\begin{lem}[$\Phi_5$]\label{lem_DOU}
The shortcut $\Phi_5$ is
$$
q_5(x;\,y;\,z)\equiv\mathrm{DOU}(x;\,y;\,z)\equiv
(\forall x_5(\mathrm{DO}(x_5;\,y;\,z)\leftrightarrow(x_5=x)))\,.
$$
Thus $|E_5|=37$, $|E_5|_{\neg}=0$ and
$$
\text{$D\models_f\mathrm{DOU}(x;\,y;\,z)$ if and only if 
$f(x)$ is the doubleton 
$\{f(y),\,f(z)\}_D$}\,.
$$
\end{lem}
\proof
In view of Lemma~\ref{lem_DO},
the syntactic part is an easy 
count $4+3+25+5=37$ and observation that no negation was used. 
Semantically, the satisfaction in $D$ of the formula 
$\mathrm{DOU}(x;y;z)$ matches the word description.
\eproof

\noindent
By the forward expansion,
{\footnotesize
$$
E_5\equiv
(\forall x_5((\forall x_4((x_4\in x_5)\leftrightarrow((x_4=y)\vee(x_4=z))))\leftrightarrow(x_5=x)))
$$
}

\begin{lem}[$\Phi_6$]\label{lem_OPA}
The shortcut $\Phi_6$ is
\begin{eqnarray*}
&&q_6(x;\,y;\,z)\equiv\mathrm{OPA}(x;\,y;\,z)\equiv
(\exists x_6(\exists x_7(\mathrm{DOU}(x;\,x_6;\,x_7)\wedge\\
&&(\mathrm{SIN}(x_6;\,y)\wedge\mathrm{DOU}(x_7;\,y;\,z)))))\,.
\end{eqnarray*}
Thus $|E_6|=117$, $|E_6|_{\neg}=0$ and 
$$
\text{$D\models_f\mathrm{OPA}(x;\,y;\,z)$ if and only if $f(x)$ is $\langle f(y),\,f(z)\rangle_D$}\,.
$$
\end{lem}
\proof
In view of Lemmas~\ref{lem_SIN} and \ref{lem_DOU},
the syntactic part is an easy 
count $4+4+3+37+3+29+37=117$ and observation that no negation was used. Semantically, the satisfaction in $D$ of the 
formula $\mathrm{OPA}(x;y;z)$ matches the word description.
\eproof

\noindent
By the forward expansion,
{\footnotesize
\begin{eqnarray*}
&&E_6\equiv(\exists x_6(\exists x_7((\forall x_5((\forall x_4((x_4\in x_5)\leftrightarrow((x_4=x_6)\vee(x_4=x_7))))\leftrightarrow(x_5=x)))
\wedge\\
&&((\forall x_3((\forall x_2((x_2\in x_3)\leftrightarrow(x_2=y)))\leftrightarrow(x_3=x_6)))\wedge\\
&&(\forall x_5((\forall x_4((x_4\in x_5)\leftrightarrow((x_4=y)\vee(x_4=z))))\leftrightarrow(x_5=x_7)))))))\,.
\end{eqnarray*}
}

Before we get to the last three shortcuts we show 
that the expansion $E_6$ of 
$\Phi_6$ along $U_0$ works as 
an ordered pair.  By Proposition~\ref{prop_freeVar}, 
$$
E_6(x,\,y,\,z)\equiv E_6
$$
is a~ZF$'$ formula with free variables $x$, $y$ and $z$.

\begin{prop}\label{prop_OPA_D}
We have
$$
D\models_f
((E_6(x,\,y,\,z)\wedge E_6(a,\,b,\,c))\to((x=a)\leftrightarrow((y=b)\wedge(z=c))))
$$
where $E_6(a,b,c)\equiv\mathrm{
sub}_1(E_6,x/a,y/b,z/c)$.
\end{prop}
\proof
This is immediate from Proposition~\ref{prop_OPA}. 
\eproof

\begin{lem}[$\Phi_7$]\label{lem_REL}
The shortcut $\Phi_7$ is
\begin{eqnarray*}
&&q_7(x;\,y)\equiv
\mathrm{REL}(x;\,y)\equiv
(\forall x_8((x_8\in x)\to(\exists x_9(\exists x_{10}(\mathrm{OPA}(x_8;\,x_9;\,x_{10})\wedge\\
&&((x_9\in y)\wedge\mathrm{SUS}(x_{10};\,y)))))))\,.
\end{eqnarray*}
Thus $|E_7|=165$, $|E_7|_{\neg}=0$ and
$$
\text{$D\models_f\mathrm{REL}(x;\,y)$ iff $f(x)$ is 
a~$D$-relation from $f(y)$ to $P(f(y))$}\,.
$$
\end{lem}
\proof
In view of Lemmas~\ref{lem_OPA} and \ref{lem_SUS},
the syntactic part is an easy 
count $4+3+5+4+4+3+117+3+5+17=165$ and observation that no negation was used. Semantically, the satisfaction in $D$ of the 
formula $\mathrm{REL}(x;y)$ matches the word description.
\eproof

\noindent
By the forward expansion,
{\footnotesize
\begin{eqnarray*}
&&E_7\equiv(\forall x_8((x_8\in x)\to(\exists x_9(\exists x_{10}\\
&&((\exists x_6(\exists x_7((\forall x_5((\forall x_4((x_4\in x_5)\leftrightarrow((x_4=x_6)\vee(x_4=x_7))))\leftrightarrow(x_5=x_8)))
\wedge\\
&&((\forall x_3((\forall x_2((x_2\in x_3)\leftrightarrow(x_2=x_9)))\leftrightarrow(x_3=x_6)))\wedge\\
&&(\forall x_5((\forall x_4((x_4\in x_5)\leftrightarrow((x_4=x_9)\vee(x_4=x_{10}))))\leftrightarrow(x_5=x_7)))))))\wedge\\
&&((x_9\in y)\wedge(\forall x_1((x_1\in x_{10})\to(x_1\in y)))))))))\,.
\end{eqnarray*}
}

\begin{lem}[$\Phi_8$]\label{lem_FUN}
The shortcut $\Phi_8$ is
\begin{eqnarray*}
&&q_8(x;\,y)\equiv
\mathrm{FUN}(x;\,y)\equiv
(\mathrm{REL}(x;\,y)\wedge(\forall x_{11}((x_{11}\in y)\to\\
&&(\exists x_{12}(\forall x_{13}((x_{13}=x_{12})\leftrightarrow((x_{13}\in x)\wedge(\exists x_{14}\mathrm{OPA}(x_{13};\,x_{11};\,x_{14})))))))))
\,.
\end{eqnarray*}
Thus $|E_8|=325$, $|E_8|_{\neg}=0$ and
$$
\text{$D\models_f\mathrm{FUN}(x;\,y)$ iff $f(x)$ is 
a~$D$-function from $f(y)$ to $P(f(y))$}\,.
$$
\end{lem}
\proof
In view of Lemmas~\ref{lem_REL} and \ref{lem_OPA},
the syntactic part is an easy 
count $3+165+4+3+5+4+4+3+5+3+5+4+117=325$ and observation that no negation was used. Semantically, the satisfaction in $D$ of the 
formula $\mathrm{FUN}(x;y)$ matches the word description.
\eproof

\noindent
By the forward expansion,
{\footnotesize
\begin{eqnarray*}
&&E_8\equiv
((\forall x_8((x_8\in x)\to(\exists x_9(\exists x_{10}\\
&&((\exists x_6(\exists x_7((\forall x_5((\forall x_4((x_4\in x_5)\leftrightarrow((x_4=x_6)\vee(x_4=x_7))))\leftrightarrow(x_5=x_8)))
\wedge\\
&&((\forall x_3((\forall x_2((x_2\in x_3)\leftrightarrow(x_2=x_9)))\leftrightarrow(x_3=x_6)))\wedge\\
&&(\forall x_5((\forall x_4((x_4\in x_5)\leftrightarrow((x_4=x_9)\vee(x_4=x_{10}))))\leftrightarrow(x_5=x_7)))))))\wedge\\
&&((x_9\in y)\wedge(\forall x_1((x_1\in x_{10})\to(x_1\in y)))))))))\wedge(\forall x_{11}((x_{11}\in y)\to\\
&&(\exists x_{12}(\forall x_{13}((x_{13}=x_{12})\leftrightarrow((x_{13}\in x)\wedge\\
&&(\exists x_{14}(\exists x_6(\exists x_7((\forall x_5((\forall x_4((x_4\in x_5)\leftrightarrow((x_4=x_6)\vee(x_4=x_7))))\leftrightarrow(x_5=x_{13})))
\wedge\\
&&((\forall x_3((\forall x_2((x_2\in x_3)\leftrightarrow(x_2=x_{11})))\leftrightarrow(x_3=x_6)))\wedge\\
&&(\forall x_5((\forall x_4((x_4\in x_5)\leftrightarrow((x_4=x_{11})\vee(x_4=x_{14}))))\leftrightarrow(x_5=x_7)))))))))))))))
\,.
\end{eqnarray*}
}

\begin{lem}[$\Phi_9$]\label{lem_SUR}
The shortcut $\Phi_9$ is
\begin{eqnarray*}
&&q_9(x;\,y)\equiv
\mathrm{SUR}(x;\,y)\equiv
(\mathrm{FUN}(x;\,y)\wedge(\forall x_{15}(\mathrm{SUS}(x_{15};\,y)\to\\
&&(\exists x_{16}(\exists x_{17}((x_{16}\in x)\wedge\mathrm{OPA}(x_{16};\,x_{17};\,x_{15})))))))\,.
\end{eqnarray*}
Thus $|E_9|=485$, $|E_9|_{\neg}=0$ and
$$
\text{$D\models_f\mathrm{SUR}(x;\,y)$ iff $f(x)$ is 
a~$D$-surjection from $f(y)$ to $P(f(y))$}\,.
$$
\end{lem}
\proof
In view of Lemmas~\ref{lem_FUN} and \ref{lem_OPA},
the syntactic part is an easy 
count $3+325+4+3+17+4+4+3+5+117=485$ and observation that no negation was used. Semantically, the satisfaction in $D$ of the 
formula $\mathrm{SUR}(x;y)$ matches the word description.
\eproof

\noindent
By the forward expansion,
{\footnotesize
\begin{eqnarray*}
&&E_9\equiv
(((\forall x_8((x_8\in x)\to(\exists x_9(\exists x_{10}\\
&&((\exists x_6(\exists x_7((\forall x_5((\forall x_4((x_4\in x_5)\leftrightarrow((x_4=x_6)\vee(x_4=x_7))))\leftrightarrow(x_5=x_8)))
\wedge\\
&&((\forall x_3((\forall x_2((x_2\in x_3)\leftrightarrow(x_2=x_9)))\leftrightarrow(x_3=x_6)))\wedge\\
&&(\forall x_5((\forall x_4((x_4\in x_5)\leftrightarrow((x_4=x_9)\vee(x_4=x_{10}))))\leftrightarrow(x_5=x_7)))))))\wedge\\
&&((x_9\in y)\wedge(\forall x_1((x_1\in x_{10})\to(x_1\in y)))))))))\wedge(\forall x_{11}((x_{11}\in y)\to\\
&&(\exists x_{12}(\forall x_{13}((x_{13}=x_{12})\leftrightarrow((x_{13}\in x)\wedge\\
&&(\exists x_{14}(\exists x_6(\exists x_7((\forall x_5((\forall x_4((x_4\in x_5)\leftrightarrow((x_4=x_6)\vee(x_4=x_7))))\leftrightarrow(x_5=x_{13})))
\wedge\\
&&((\forall x_3((\forall x_2((x_2\in x_3)\leftrightarrow(x_2=x_{11})))\leftrightarrow(x_3=x_6)))\wedge\\
&&(\forall x_5((\forall x_4((x_4\in x_5)\leftrightarrow((x_4=x_{11})\vee(x_4=x_{14}))))\leftrightarrow(x_5=x_7)))))))))))))))\wedge\\
&&(\forall x_{15}((\forall x_1((x_1\in x_{15})\to(x_1\in y)))\to\\
&&(\exists x_{16}(\exists x_{17}((x_{16}\in x)\wedge\\
&&(\exists x_6(\exists x_7((\forall x_5((\forall x_4((x_4\in x_5)\leftrightarrow((x_4=x_6)\vee(x_4=x_7))))\leftrightarrow(x_5=x_{16})))
\wedge\\
&&((\forall x_3((\forall x_2((x_2\in x_3)\leftrightarrow(x_2=x_{17})))\leftrightarrow(x_3=x_6)))\wedge\\
&&(\forall x_5((\forall x_4((x_4\in x_5)\leftrightarrow((x_4=x_{17})\vee(x_4=x_{15}))))\leftrightarrow(x_5=x_7)))))))))))))\,.
\end{eqnarray*}
}

By performing the previous expansions we actually proved that 
the defined abbreviation scheme is well formed, but let us recapitulate
it.

\begin{prop}\label{prop_regAbbr}
The abbreviation scheme 
$$
U_0\equiv\langle\Phi_1,\,\Phi_2,\,\ds,\,\Phi_9\rangle
$$ 
defined by Lemmas~\ref{lem_SUS}--\ref{lem_OPA} and \ref{lem_REL}--\ref{lem_SUR} is well formed.    
\end{prop}
\proof
Indeed, $R(\Phi_1)=
\emptyset$, $R(\Phi_2)=
\emptyset$, $R(\Phi_3)=\{2\}$,
$R(\Phi_4)=\emptyset$, $R(\Phi_5)=
\{4\}$, $R(\Phi_6)=
\{3,5\}$, $R(\Phi_7)=\{1,6\}$,
$R(\Phi_8)=\{6,7\}$, $R(\Phi_9)=
\{1,6,8\}$, $V(\Phi_1)=\{1\}$, $V(\Phi_2)=
\{2\}$, $V(\Phi_3)=\{3\}$, 
$V(\Phi_4)=\{4\}$, 
$V(\Phi_5)=\{5\}$.
$V(\Phi_6)=\{6,7\}$, 
$V(\Phi_7)=\{8,9,10\}$,
$V(\Phi_8)=\{11,12,13,14\}$, and 
$V(\Phi_9)=\{15,16,17\}$.
\eproof

We arrive at our main result.

\begin{thm}\label{thm_CantorForm}
Let
$$
\varphi\equiv
(\forall x_{18}\neg(\exists x_{19} E_9(x_{19};\,x_{18})))\ \ (\in\mathcal{A}^*)
$$
where $E_9$ is the last expansion along the abbreviation scheme $U_0$ and 
$$
E_9(x_{19};\,x_{18})\equiv
\mathrm{sub}_1(E_9,\,x/x_{19},\,y/x_{18})\,.
$$
Then
$$
|\varphi|=4+1+4+485=494,\ |\varphi|_{\neg}=1
$$
and
$$
\text{$D\models\varphi$ if and only if Cantor's theorem holds in $D$}\,.
$$
\end{thm}
\proof 
By Proposition~\ref{prop_freeVar}, $\varphi$ is a~ZF sentence. The theorem is 
an immediate corollary of Lemma~\ref{lem_SUR}. 
\eproof

\noindent
Explicitly, $\varphi$ is
{\footnotesize
\begin{eqnarray*}
&&\varphi\equiv(\forall x_{18}\neg(\exists x_{19} (((\forall x_8((x_8\in x_{19})\to(\exists x_9(\exists x_{10}\\
&&((\exists x_6(\exists x_7((\forall x_5((\forall x_4((x_4\in x_5)\leftrightarrow((x_4=x_6)\vee(x_4=x_7))))\leftrightarrow(x_5=x_8)))
\wedge\\
&&((\forall x_3((\forall x_2((x_2\in x_3)\leftrightarrow(x_2=x_9)))\leftrightarrow(x_3=x_6)))\wedge\\
&&(\forall x_5((\forall x_4((x_4\in x_5)\leftrightarrow((x_4=x_9)\vee(x_4=x_{10}))))\leftrightarrow(x_5=x_7)))))))\wedge\\
&&((x_9\in x_{18})\wedge(\forall x_1((x_1\in x_{10})\to(x_1\in y)))))))))\wedge(\forall x_{11}((x_{11}\in x_{18})\to\\
&&(\exists x_{12}(\forall x_{13}((x_{13}=x_{12})\leftrightarrow((x_{13}\in x)\wedge\\
&&(\exists x_{14}(\exists x_6(\exists x_7((\forall x_5((\forall x_4((x_4\in x_5)\leftrightarrow((x_4=x_6)\vee(x_4=x_7))))\leftrightarrow(x_5=x_{13})))
\wedge\\
&&((\forall x_3((\forall x_2((x_2\in x_3)\leftrightarrow(x_2=x_{11})))\leftrightarrow(x_3=x_6)))\wedge\\
&&(\forall x_5((\forall x_4((x_4\in x_5)\leftrightarrow((x_4=x_{11})\vee(x_4=x_{14}))))\leftrightarrow(x_5=x_7)))))))))))))))\wedge\\
&&(\forall x_{15}((\forall x_1((x_1\in x_{15})\to(x_1\in x_{18})))\to\\
&&(\exists x_{16}(\exists x_{17}((x_{16}\in x_{19})\wedge\\
&&(\exists x_6(\exists x_7((\forall x_5((\forall x_4((x_4\in x_5)\leftrightarrow((x_4=x_6)\vee(x_4=x_7))))\leftrightarrow(x_5=x_{16})))
\wedge\\
&&((\forall x_3((\forall x_2((x_2\in x_3)\leftrightarrow(x_2=x_{17})))\leftrightarrow(x_3=x_6)))\wedge\\
&&(\forall x_5((\forall x_4((x_4\in x_5)\leftrightarrow((x_4=x_{17})\vee(x_4=x_{15}))))\leftrightarrow(x_5=x_7)))))))))))))))\,.    
\end{eqnarray*}
}

\section{Extensive and strongly extensive digraphs}\label{sec_digr}

We define these two families of digraphs. Let $D=\langle V,E\rangle$ 
be a~digraph.

\begin{defi}[extensive digraphs]\label{def_extSM}
$D$ is extensive if the axiom schema of specification holds in it\,---\,for every 
$n+3$, $n\in\N$, mutually distinct variables $v_1$, $v_2$, $\ds$, $v_{n+3}$ in 
$\{x_1,x_2,\ds\}$ and every 
{\em ZF} formula $\psi$ with free occurrences only of (some of) the variables $v_1$, $v_2$, $\ds$, $v_{n+2}$,
\begin{eqnarray*}
&&D\models(\forall v_1(\ds(\forall v_n(\forall v_{n+1}(\exists v_{n+3}(\forall v_{n+2}((v_{n+2}\in v_{n+3})\leftrightarrow\\
&&((v_{n+2}\in v_{n+1})\wedge\psi))))))\ds))\,.
\end{eqnarray*}
\end{defi}
 
If $v\in V$, $u$ is a~$D$-surjection from $v$ to 
$P(v)$ and $wEv$, we write
$$
w\stackrel{u}{\mapsto}w'
$$
to denote the unique vertex $w'\sus_D v$ such that
$\langle w,w'\rangle_D$ is a~$D$-element of $u$. We obtain the 
following generalization of Cantor's Theorem~\ref{thm_Cantor}.

\begin{thm}\label{thm_MCantor}
Cantor's theorem holds in every extensive digraph.
\end{thm}
\proof
For the contrary, let $D=\langle V,E\rangle$ be an 
extensive digraph that is not Cantor: there exists a~vertex $u\in V$ and 
a~$D$-surjection $v\in V$ from $u$ to $P(u)$. We consider the set of vertices
$$
A\equiv\{w\in_D u\cc\;w\not\in N(w')\text{ for the vertex $w'$ given by }w\stackrel{v}{\mapsto}w'\}\,.
$$
$A$ is defined by an axiom of specification and $D$ is extensive, 
hence there is a~vertex $w_0\in V$ such that
$$
N(w_0)=A\,.
$$
Since $w_0\sus_D u$, there is a~vertex $w_1\in_D u$ such that 
$
w_1\stackrel{v}{\mapsto}w_0$. 
We get the contradiction that $w_1\in A$ iff $w_1\not\in N(w_0)$ iff $w_1\not\in A$.
\eproof

\noindent
A~different Cantor's theorem for digraphs was found by 
Fajtlowicz \cite{fajt}.

Models of ZF set theory are extensive and therefore Cantor.

\begin{prop}\label{prop_onDZF}
Every model $D_{\mathrm{ZF}}=\langle V,E\rangle$ of {\em ZF} set theory is extensive.     
\end{prop}
\proof
The axiom schema of specifications is an axiom schema of ZF.
\eproof
\vspace{-3mm}
\begin{defi}[strongly extensive digraphs]\label{def_strExt}
We say that a~digraph $D=\langle V,E\rangle$ is strongly extensive if 
for every vertex $u\in V$ and every 
set $A$ of in-neighbors of $u$ there exists a~vertex $v\in V$ such that $N(v)=A$.    
\end{defi}
Every strongly extensive digraph is extensive and it is easy to see that every 
finite extensive digraph is strongly extensive.

We give four examples of strongly extensive digraphs. The first three 
are finite:
$\langle[1],\,\emptyset\rangle$, 
$\langle[2],\,\{\langle 1,1\rangle\}\rangle$
and
$$
\langle[4],\,\{\langle 1,1\rangle,\,\langle 2,1\rangle,\,\langle 1,3\rangle,\,\langle 2,4\rangle\}\rangle\,.
$$ 
The fourth example is infinite and countable.

\begin{prop}\label{prop_couStrExt}
There exists a~digraph 
$D=\langle V,E\rangle$ with countable vertex set~$V$ such 
that $N(u)$ is finite for every $u\in V$ and that 
for every finite set $A\sus V$ there exists a~$u\in V$ with $N(u)=A$. In particular, $D$ is strongly extensive.    
\end{prop}
\proof
For $n\in\N$ we define finite sets $V_n\sus\N$, $V_1<V_2<\ds$, and $E_n\sus\N^2$. We start with  
$V_1\equiv\{1\}$ and $E_1\equiv\emptyset$. Suppose that 
$V_1$, $\ds$, $V_n$ and $E_1$, $\ds$, $E_n$ are already defined. We
set $X\equiv V_1\cup\ds\cup V_n$ and 
$m\equiv|\mathcal{P}(X)|$. We take any enumeration of subsets of $X$,
$$
\mathcal{P}(X)=
\{A_1,\,A_2,\,\ds,\,A_m\}\,,
$$
and set
$$
V_{n+1}\equiv[\max(V_n)+1,\,\max(V_n)+m]
$$
and
$$
E_{n+1}\equiv\bigcup_{i=1}^m
A_i\times\{\max(V_n)+i\}\,.
$$
We finally define 
$${\textstyle
D=\langle V,\,E\rangle\equiv\big\langle
\bigcup_{n=1}^{\infty}V_n,\,\bigcup_{n=1}^{\infty}E_n\big\rangle\,.
}
$$
In fact, $V=\N$. It is not hard to see that $D$ has both stated 
properties. 
\eproof

\noindent
By Theorem~\ref{thm_MCantor}, in each of the four previous digraphs
Cantor's theorem holds.

\section{Concluding remarks}\label{sec_concl}

For the previous version of this article, rich in quotations, see 
\cite{klaz_arXiv}. Let $\mathcal{F}(\mathcal{A})$ ($\sus\mathcal{A}^*$) 
be the set of 
ZF formulas and $\mathcal{CAN}$ ($\in\N$) be the minimum length $|\psi|$ of 
a~sentence $\psi\in\mathcal{F}(\mathcal{A})$ such that for every 
digraph $D$,
$$
D\models\psi\iff\text{Cantor's theorem holds in $D$}\,.
$$

\begin{prob}
Give good upper and lower bounds on $\mathcal{CAN}$ or determine this number exactly.    
\end{prob}
By Theorem~\ref{thm_CantorForm}, $\mathcal{CAN}\le 494$. Good lower bounds would be interesting. Progress
on this problem might be achieved by
obtaining some simple structural characterization of Cantor digraphs.

Let $n\in\N$, $e_n$ be the number of (finite) strongly extensive 
digraphs $D=\langle[n],E\rangle$ and $c_n$ be the number of digraphs 
$D=\langle[n],E\rangle$ 
such that Cantor's theorem holds in $D$.

\begin{prob}
Give good upper and lower bounds on $e_n$ and $c_n$ or determine these numbers exactly.    
\end{prob}
By Theorem~\ref{thm_MCantor}, $e_n\le c_n$. How much 
larger than $e_n$ is $c_n$? Efficient characterizations of both kinds of digraphs would be interesting.  

Besides Cantor's theorem there is the Cantor--Bernstein theorem: for 
any sets $x$
and $y$, if there exists an injection from $x$ to $y$, and an injection from $y$ to 
$x$, then there exists a~bijection from $x$ to $y$. It is possible to 
do to the C.--B. theorem what we did to Cantor's theorem. What are other interesting theorems on sets and 
functions between them? 

G\"odel's completeness theorem (GCT) \cite{gode}, see also \cite{henk} 
and \cite[Chapter~IV]{daws}, says that
$$
\text{a~FO theory $\mathcal{T}$ is consistent $\iff$ $\mathcal{T}$ has a~model $\mathcal{M}$.}    
$$
That is, one cannot deduce a~contradiction from
the set $\mathcal{T}$ of sentences stated in a~first order language if 
and only if there exists a~structure 
$\mathcal{M}$ such that every sentence $\psi$ in $\mathcal{T}$ is 
true in $\mathcal{M}$, i.e. $\mathcal{M}\models\psi$. We view 
the right-hand side of the equivalence as too 
informal. Sentences in $\mathcal{T}$ are very precise objects, namely certain 
words, but $\mathcal{M}$ is only  a~naive set universe with naive sets of tuples. This article started as 
a~project aiming at obtaining more 
rigorous statements of GCT, and Cantor's 
theorem was to be only an illustration. Now we understand the 
nature of our GCT project better. The statement of Cantor's 
theorem by sentence $\varphi$ in Theorem~\ref{thm_CantorForm} is, of 
course, still ``model-vague'', it involves digraphs $D$, but there cannot be anything more rigorous and 
precise than the word $\varphi$. 
In \cite{klaz} we hope to obtain an analogous sentence for GCT. 
This will be obviously much harder than what we did for Cantor's theorem, but we are 
confident that it can be done.
Why? Because it has been already done, only in a~different language than we use here\,---\,the proof of GCT and its statement were formalized by 
From \cite{from}.

\appendix

\section{Zkratky}\label{app_zkr}

I went from PLR to MLR through 
\v{C}SSR, SNB stopped my DKW for TK.

\bigskip\noindent
PLR is an acronym for {\em Polsk\'a lidov\'a republika} or {\em the 
People's Republic of Poland}, which was the official Czech name for 
Poland in 1947--1989. MLR is an acronym for {\em Mad'arsk\'a 
lidov\'a republika} or {\em the People's Republic of Hungary}, which was the official 
Czech name for Hungary in 1949--1989.
\v{C}SSR is an acronym for 
{\em \v{C}eskoslovensk\'a socialistick\'a republika} or {\em 
the Czechoslovak Socialistic Republic}, which was in 1960--1990 the official 
Czech name for the state composed of the present Czechia and Slovakia. In 1948--1960 the official name was just 
{\em \v{C}eskoslovensk\'a republika}, with the acronym \v CSR. You see, we were more 
advanced than Poland or Hungary. SNB  is an acronym for {\em Sbor 
n\'arodn\'\i\ bezpe\v cnosti}. This is not worth translating to 
English, it was communist police. 
DKW is an acronym for the German
{\em Dampfkraftwagen} or {\em steam car}\,---\,see \cite{dkw} for more information. Finally, TK is an acronym for {\em 
techick\'a kontrola} or {\em 
technical check}.

\bigskip\noindent
{\em Department of Applied Mathematics\\
Faculty of Mathematics and Physics\\
Charles University\\
Malostransk\'e n\'am. 25\\
118 00 Praha\\
Czechia}\\
{\tt klazar@kam.mff.cuni.cz}

\end{document}